\documentclass[10pt]{amsart}
\usepackage{amsmath,graphicx,amscd}

\theoremstyle{plain}
\newtheorem{theorem}{Theorem}[section]

\newtheorem{corollary}[theorem]{Corollary}
\newtheorem{proposition}[theorem]{Proposition}

\theoremstyle{remark}
\newtheorem*{remark}{Remark}

\theoremstyle{definition}

\def \R {\mathbf{R}}

\def \E {\mathbb{E}}
\def \S {\mathbb{S}}

\def\PP{\mathcal{P}}

\def\RR{\mathcal{R}}

\def\MM{\mathcal{M}}
\numberwithin{equation}{section}
\def\SOf{{\rm SO}(4)}
\def\SO{{\rm SO}(3)}
\def\SU{{\rm SU}(2)}

\def\id{{\rm Id}}

\def\fg{\pi_1}

\begin{document}
\title{On the symplectic volume  of the moduli space of Spherical and Euclidean polygons}
\author {Vu The Khoi}
\address{ Institute of Mathematics, 18 Hoang Quoc Viet road, 10307, Hanoi, Vietnam}
\email{vtkhoi@math.ac.vn}
\thanks{The author was partially supported by the National Basic Research Program of Vietnam}
\subjclass{  53D30.}
\keywords{moduli space of polygons, symplectic volume}

\begin{abstract}
In this paper, we study  the symplectic volume of the moduli space of polygons   by using Witten's formula.
We propose to use this volume as a measure for the flexibility of a polygon with fixed side-lengths. 
 The main result of our is  that among all the polygons with fixed perimeter in $\S^3$ or $\E^3$ the regular one is the most flexible and that among all the spherical polygons   the  regular one with side-length $\pi/2$ is the most flexible.      
\end{abstract}
\maketitle
\section{Introduction} 
\vskip0.1cm
An polygon in $\S^3$ or $\E^3$ is specified by its set of vertices $v=(v_1,\ldots,v_n).$ This vertices are joined in cyclic order by edges $e_1,\ldots,e_n$, where $e_i$ is the directed geodesic segment from $v_i$ to $v_{i+1}.$ Two polygons   $P=(v_1,\ldots,v_n)$ and $Q=(w_1,\ldots,w_n)$ are identified if there exists an orientation preserving isometry sending each $v_i$ to $w_i.$ 

The side-length $r_i$ of a polygon is defined to be the length of the geodesic segment $e_i.$
We say that a polygon is \textit {regular} if all of its side-lengths are equal.

Let $r=(r_1,\ldots,r_n)$ be a tuple of real numbers such that $0<r_i<\pi\  \forall i,$ following \cite{mp},  we will denote by   $\PP_r^{\S^3}$ the  \textit {configuration space of all polygons   in $\S^3$ with side-lengths} $r.$   The \textit {moduli space of polygons   in $\S^3$ with side-lengths} $r$ is defined to be $\PP_r^{\S^3}/\mathrm{Iso}=\PP_r^{\S^3}/\SOf$. 

 In $\E^3,$ similarly,  if $r=(r_1,\ldots,r_n)$ is a tuple of positive real numbers, we denote by  $\PP_r^{\E^3}$ the  \textit {configuration space of all  polygons   in $\E^3$ with side-lengths} $r.$ The \textit {moduli space of polygons   in $\E^3$ with side-lengths} $r$ is defined to be $\PP_r^{\E^3}/\mathrm{Iso}=\PP_r^{\E^3}/\E^+(3).$ 
  
We will drop the superscript and simply write $\PP_r$ when we want to talk about polygons in both $\S^3$ and $\E^3.$  
  
We are interested in the question : \textit {for which tuple $r$ the polygon with side-lengths $r$ is the most flexible when it is   moved  in the space provided that its  side-lengths  are fixed} ?  First of all we need to specified what is the exact meaning of  ``flexible". When a polygon is moved, it always lies in its configuration space. Therefore, it is natural to think that a polygon with side-lengths $r$ is more flexible than another one with side-lengths $r'$ if the corresponding configuration space $\PP_r$ is ``bigger" than $\PP_{r'}$ in some sense. Fortunately, there is a standard symplectic structure on the moduli spaces $\PP_r/\mathrm{Iso}$, see \cite{ab,g}, and \textit {we will measure the flexibility  of a polygon with side-lengths $r$ by the symplectic volume of its moduli space}.

The main results are the followings:
\begin{theorem} Among all the  spherical or Euclidean polygons with fixed perimeter, the regular polygon is one of the most flexible. Moreover in the case $n$ is even, the regular $n$-gon is the unique one with this property. 
\end{theorem}
\begin{theorem} Among all the  spherical polygons the regular one with side-length $\pi/2$ is the unique one which is the most flexible.  
\end{theorem}
The rest of this paper is organized as follows. In section 2, we will recall the symplectic structure on the moduli space $\PP_r^{\S^3}/\SOf$ and using Witten's formula to find its symplectic volume in terms of Bernoulli polynomials.  As a corollary we also derive the positivity of certain trigonometric sums which has been proved by analytic method.
Section 3 is devoted to the study of the moduli space $\PP_r^{\E^3}/\E^+(3),$ in particular we derive an explicit formula for its symplectic volume.
 The main results are proved in section 4 by a detail study of the symplectic volume of as a function of $r.$  
\section{ Volume of the moduli space of spherical polygons  } 
\vskip0.1cm   
We first briefly recall the identification between $\PP_r^{\S^3}/\SOf$ and the moduli space of flat $\SU$ connections on a punctured sphere with fixed holonomies around the punctures. Let $S_n:=S^2\setminus \{p_1, \ldots, p_n\}$ be an $n$-punctured sphere. For a tuple of numbers $r=(r_1,\ldots ,r_n)$  such that $0<r_i<\pi\  \forall i,$ we will denote by $\MM(S_n,r)$  the moduli space of flat $\SU$ connections on $S_n$ modulo gauge equivalence such that the holonomy around $p_i$ is conjugate  to   $\left(\begin{array}{cc}e^{ir_i}&0\\0&e^{-ir_i}\end{array}\right).$

It is well-known that $\MM(S_n,r)\cong \RR(S_n,r)/\SU,$ where
$$\RR(S_n,r):=\{(g_1,\ldots, g_n)\in \SU^n| \ g_1\cdots g_n=\mathrm {Id}, \mathrm{Tr}(g_j)=2\cos(r_j)\  \forall j\}$$
   is the representation space consists of all representations from $\fg(S_n)$ to $\SU$ such that the image of the loop around the puncture $p_i$ is conjugate to $\left(\begin{array}{cc}e^{ir_i}&0\\0&e^{-ir_i}\end{array}\right)$
and  the group $\SU$  act  diagonally by  conjugation.

On the other hand,  following \cite{mp}, we may identify $\RR(S_n,r)$ with the configuration space $\PP_{r,0}^{\S^3}$ of  \textit {based polygons,} i.e. polygons having the first vertex $v_1= \id\in \SU \cong \S^3.$ For $(g_1,\ldots, g_n) \in \RR(S_n,r)$ we construct a based polygon by defining $g_0 := \id$ and let $v_j:=g_0 \cdots g_{j-1}$  for $1\le j \le n.$ Under this identification, we have $$\RR(S_n,r)/\SU\cong \PP_{r,0}^{\S^3}/\SO \cong \PP_{r}^{\S^3}/\SOf.$$ 
Therefore, we have a natural identification  $\PP_{r}^{\S^3}/\SOf \cong \MM(S_n,r).$ According to Atiyah-Bott \cite{ab}, the moduli space  of flat  connections on a surface of genus $g$ with specified holonomies around the boundary is a finite dimensional symplectic manifold, possibly with singularities. Witten \cite{w} compute the symplectic volume of the  moduli space of flat connections. Since then, other people have given simplified and rigorous proofs (see \cite{d,jw,liu, mw}).      
As a special case of Witten's theorem, we have the followings (see \cite{w, mw}):
\begin{theorem}(Witten's formula)  Denote by $V^{\S^3}(r)$ the symplectic volume of $\PP_{r}^{\S^3}/\SOf $- the moduli space of polygons of side-lengths $r=(r_1,\ldots,r_n).$ 
 We have:
$$V^{\S^3}(r) = \frac{2^{n-1}}{\pi}\sum_{k=1}^\infty \frac{\sin kr_1 \ldots \sin kr_n}{k^{n-2}}. $$ 
\end{theorem}

For $I\subset \{1,2, \ldots n \},$ let $\bar I$ denote the complement of $I, $ $| I |$ be the cardinality of $I$ and $r_I:=\sum_{i\in I} r_i.$ It has been shown that the moduli space of polygons of side-lengths $r=(r_1,\ldots,r_n)$ is non-empty if and only if :
$$r_I \le r_{\bar I} + (|I| - 1)\pi \quad \forall I\subset \{1,2, \ldots n \}, |I| \quad \textrm{odd}.$$ 

See \cite{aw} for a proof using quantum cohomology. This result can also be proved in a similar line as in the case of spherical polygons in $\S^2$ in \cite{km}.

Combining this result with Witten's formula we obtained the positivity of a trigonometric series.
\begin{corollary} If $r=(r_1,\ldots,r_n)$ be a tuple of real numbers such that $0<r_i<\pi $   for all $i,$ Then 
$$ \sum_{k=1}^\infty \frac{\sin kr_1 \ldots \sin kr_n}{k^{n-2}}\ge 0 .$$
Moreover, the inequality is strict  if and only if $r_I \le r_{\bar I} + (|I| - 1)\pi$  for all $I\subset \{1,2, \ldots n \}$, $|I|$  odd. 
\end{corollary}
      
We can find a closed form expression of the volume function in terms of the Bernoulli polynomials.
\begin{proposition} \label{prop1}
The symplectic volumes of the moduli space of  $2n$-gons and $(2n+1)$-gons  are given respectively by:
$$V^{\S^3}(r) = \frac{(2\pi)^{2n-3}}{2(2n-2) !}\sum_{I\subset \{1,2, \ldots 2n \}}(-1)^{| I |}B_{2n-2}(\big \{(\sum_{i\in I} r_i-\sum_{i\in \bar I} r_i)/2\pi\big \})$$ 
$$V^{\S^3}(r) = \frac{(2\pi)^{2n-2}}{(2n-1) !}\sum_{I\subset \{1,2, \ldots 2n+1 \},  |I|\  \textrm{odd}}B_{2n-1}(\big \{(\sum_{i\in I} r_i-\sum_{i\in \bar I} r_i)/2\pi\big \}).$$ 
Where, $B_n$ are the Bernoulli polynomials defined by the generating functions
$$\frac{te^{xt}}{e^t-1} = \sum_{n=0}^\infty B_n(x)\frac{t^n}{n !}$$
  and $\{\ \}$ denote  the fractional part.  
\end{proposition}   
\begin{proof}
In the proof of this proposition we will use the following formulae ( see \cite{gr}, 1.443) :
$$\sum_{k=1}^\infty \frac{\cos k\pi x}{k^{2n}}= (-1)^{n-1}\frac{(2\pi)^{2n}}{2(2n)!}B_{2n}(\frac{x}{2}) \qquad (1)$$
 $$\sum_{k=1}^\infty \frac{\sin k\pi x}{k^{2n+1}}= (-1)^{n-1}\frac{(2\pi)^{2n+1}}{2(2n+1)!}B_{2n+1}(\frac{x}{2}) \qquad (2)$$
Where, $B_n$ is the Bernoulli polynomials defined above and  the formulae hold for 
 $0\le x \le 2. $ 

We first prove the $2n$-gon case. Notice that we can write,
\begin{eqnarray*}
\sin kr_1 \ldots \sin kr_{2n} &= &\frac{(e^{ikr_1} - e^{-ikr_1}) \ldots (e^{ikr_{2n}} - e^{-ikr_{2n}}) } {(2i)^{2n}} \\
  &=& \frac{\sum_{I\subset \{1,2, \ldots 2n \}}  (-1)^{|I'|}e^{ik(\sum_{i\in I}r_i - \sum_{i\in I'}r_i)}} {(2i)^{2n}}\\
&=&  \frac{\sum_{I\subset \{1,2, \ldots 2n \} }  (-1)^{|I|} \cos k(\sum_{i\in I}r_i - \sum_{i\in I'}r_i)} {(2i)^{2n}}.
\end{eqnarray*}
In the last line, we have replaced $|I'|$ by $|I|$ since $|I| + |I'|=2n.$ 
Substitute this in to Witten's formula we get:
$$V^{\S^3}(r) = \frac{(-1)^n} {2\pi}\sum_{I\subset \{1,2, \ldots 2n \}}  (-1)^{|I|}\sum_{k=1}^\infty \frac{  \cos k(\sum_{i\in I}r_i - \sum_{i\in I'}r_i) }  {k^{2n-2}}$$    
Now using (1) above, we can deduce the result.

Now consider the $2n+1$-gon case. Similarly, we can write 
\begin{eqnarray*}
\sin kr_1 \ldots \sin kr_{2n+1} &= &\frac{(e^{ikr_1} - e^{-ikr_1}) \ldots (e^{ikr_{2n+1}} - e^{-ikr_{2n+1}}) } {(2i)^{2n+1}} \\
  &=& \frac{\sum_{I\subset \{1,2, \ldots 2n+1 \}}  (-1)^{|I'|}e^{ik(\sum_{i\in I}r_i - \sum_{i\in I'}r_i)}} {(2i)^{2n+1}}\\
&=&  \frac{\sum_{I\subset \{1,2, \ldots 2n+1 \}, |I|\  \textrm{odd}}   \sin k(\sum_{i\in I}r_i - \sum_{i\in I'}r_i)} {(2i)^{2n}}.
\end{eqnarray*}
 We   can obtain the formula by plugging this into Witten's formula and  then use formula (2).  

\end{proof}
\section{ Volume of the moduli space of Euclidean polygons  } 
\vskip0.1cm   
In Euclidean space, the moduli space of polygons   $\PP_r^{\E^3}/\E^+(3)$ also has a symplectic structure \cite{km, mp} which can be briefly described as follows.

Let $r=(r_1,\ldots,r_n)$ be a tuple of positive numbers, denote by  $\PP_{r,0}^{\E^3}$ the configuration space of  \textit {based polygons,} i.e. polygons having the first vertex $v_1 = 0\in \E^3.$ Consider the map:
\begin{eqnarray*}
\Phi_r:\ \PP_{r,0}^{\E^3}\qquad &\longrightarrow &\qquad (S^2)^n \\
  (v_1,\ldots,v_n) &\mapsto & (\frac{v_2}{r_1}, \frac{v_3-v_2}{r_2}, \ldots, \frac{v_n-v_{n-1}}{r_{n-1}}, \frac{-v_n}{r_n}). \\ 
\end{eqnarray*}   
If we defined $\mu_r : (S^2)^n \rightarrow \R : (u_1,\ldots,u_n) \mapsto r_1u_1 + \cdots + r_nu_n$, then $\Phi_r$ define an $\SO$-equivariant map between $\PP_{r,0}^{\E^3}$ and $\mu_r^{-1}(0)$ which results in a diffeomorphism 
  between $\PP_{r,0}^{\E^3}/\SO$  and $ \mu_r^{-1}(0)/\SO.$      

  It is well-known that 
the unit sphere  $S^2$ is a symplectic manifold where the symplectic $2$-form is the volume form. 
Consider the symplectic manifold  $(S^2)^n$ with the weighted symplectic structure $r_1, \ldots, r_n$, i.e., the symplectic $2$-form is  $r_1\omega_1 + \cdots + r_n\omega_n, $ where $\omega_i$ is the volume form on the $i^{th}$ component. 
It is well-known that $\mu_r$ is the moment map for $(S^2)^n$ with the diagonal action of $\SO.$  Therefore, we can identify $\PP_r^{\E^3}/\E^+(3) \cong  \PP_{r,0}^{\E^3}/\SO$ with the weighted symplectic quotient  of  $(S^2)^n$ by the diagonal action of $\SO.$ 

Moreover, we can deduce from a result of L. Jeffrey \cite{lisa}, Theorem 6.6 (see also the discussion at the end of \cite{mp})  that the moduli space of polygons    in $\S^3$ and $\E^3$ with the same  side-lengths $r$ are symplectomorphic provided that $r_i$ are sufficiently small. So the volume of the moduli space of polygons   in $\S^3$ and $\E^3$ are the same if the side-lengths are sufficiently small. 
Using this fact,  we can express   the volume  of  $\PP_r^{\E^3}/\E^+(3)$ in a quite explicit form.
 \begin{proposition}\label{prop2}
The symplectic volumes of the moduli space of  normalized Euclidean $2n$-gon and $(2n+1)$-gon  are given respectively by:
$$V^{\E^3}(r) = -\frac{1}{4(2n-3) !}\sum_{I\subset \{1,2, \ldots 2n \}}(-1)^{| I |} |\sum_{i\in I} r_i-\sum_{i\in \bar I} r_i |^{2n-3}$$ 
$$V^{\E^3}(r) = -\frac{1}{2(2n-2) !}\sum_{I\subset \{1,2, \ldots 2n+1 \},  |I|\  \textrm{odd}}\mathrm{sign}(\sum_{i\in I} r_i-\sum_{i\in \bar I} r_i) (\sum_{i\in I} r_i-\sum_{i\in \bar I} r_i)^{2n-2}.$$ 

\end{proposition}   
\begin{proof}
Notice that if we scale $r$ by a  scalar  $\lambda >0$ then  $\PP_{\lambda r}^{\E^3}/\E^+(3)$ is the weighted symplectic quotient of $(S^2)^n$ for the weight $\lambda r.$ Moreover, we know that   $\PP_r^{\E^3}/\E^+(3)$ is of dimension $2(n-3)$, it  follows that $V^{\E^3}(\lambda r) = \lambda^{n-3}V^{\E^3}(r).$ 

We first, assume that all $r_i$ are sufficiently small, from the discussion above, we know that the volume of  $\PP_r^{\E^3}/\E^+(3)$ can also be expressed in terms of Bernoulli polynomials as in Proposition \ref{prop1}.
When  all $r_i$ are  are sufficiently small, we have:
$$\{(\sum_{i\in I} r_i - \sum_{i\in \bar I} r_i)/2\pi \}= 
 {\Bigg \{} \begin{matrix} 
 (\sum_{i\in I} r_i - \sum_{i\in \bar I} r_i)/2\pi,  &  {\textrm if}\   \sum_{i\in I} r_i - \sum_{i\in \bar I} r_i \ge 0 \\
               1 + (\sum_{i\in I} r_i - \sum_{i\in \bar I} r_i)/2\pi,  &  {\textrm if}\  \sum_{i\in I} r_i - \sum_{i\in \bar I} r_i < 0. \\ \end{matrix}$$
  We recall the  two following property of Bernoulli polynomials \cite{we}:

a) $B_n(1-x) = (-1)^nB_n(x) $

b) $B_n(x)= (B+x)^n$, where the symbol $B^n$  is the Bernoulli number  $B_n(0).$ 

By the homogeneity of $V^{\E^3}(r)$, we know that to find $V^{\E^3}(r)$ it is enough to collect the terms of order 
$(n-3)$ from the formulae of Proposition \ref{prop1}.  Now using this two properties of  Bernoulli polynomials and notice that $B_1=-1/2$, one can easily obtain the required results. So our formulae hold  in the case where all $r_i$ are sufficiently small. But as notice before, $V^{\E^3}(r)$ are homogenous of degree $n-3$ with respect to multiplication by a positive scalar, therefore the formulae hold for all $r.$
\end{proof}
As a special case of our result, when all $r_i=1,$ we obtain the following formula which has been shown by direct computation in \cite{kam}. 
\begin{corollary} The volume of the moduli space of regular polygon with side-length $1$ in $\E^3$ is given by:
$$V^{\E^3}(1) = -\frac{1}{2(n-3) !} \sum_{k=0}^{[n/2]} (-1)^k \binom {n} {k}  (n-2k)^{n-3}.$$      
\end{corollary}   
\begin{proof}
 The $2n$-gon case is straightforwards, we will give the proof for the case of $2n+1$-gon.
From Proposition \ref{prop2}, we get 
 $$V^{\E^3}(1) = -\frac{1}{2(2n-2) !} ( \sum_{k \ {\textrm odd}, k > n}  \binom {2n+1} {k}  (2n+1-2k)^{2n-2}- \sum_{k\ {\textrm odd}, k\le n}  \binom {2n+1} {k}  (2n+1-2k)^{2n-2}).$$
By changing the index in the first sum from $k$ to $(2n+1-k)$ we get,  

$$V^{\E^3}(1)  = -\frac{1}{2(2n-2) !} ( \sum_{k \ {\textrm even}, k \le n}  \binom {2n+1} {k}  (2n+1-2k)^{2n-2}- \sum_{k\ {\textrm odd}, k\le n}  \binom {2n+1} {k}  (2n+1-2k)^{2n-2})$$
$ \qquad \quad  \ = -\frac{1}{2(2n-2) !}  \sum_{k=0} ^{n}  (-1)^k \binom {2n+1} {k}  (2n+1-2k)^{2n-2}.$

This finishes the proof.
\end{proof}
\section{Proof of the main results}
In this section we will give the proofs of Theorem 1.1 and 1.2 by study the maximum of the volume function which is the trigonometric series given by Theorem 2.1. From Proposition \ref{prop1}, we can easily find out that for $n=3$, the moduli space of polygons has volume $1,$ no matter what the side-lengths are. So  we only need to deal with the non-trivial cases and may assume that $n>3.$   

\textit {Proof of Theorem 1.1.}   
It is enough to prove the spherical case.

For a fixed $n,$ consider the set of all polygon whose perimeter is $P.$ To prove the first part of our theorem, we must show the following :

\textit{Claim:} Let  $$f(x_1, \ldots, x_n) = \sum_{k=1}^\infty \frac{\sin kx_1 \ldots \sin kx_n}{k^{n-2}}$$ be a function 
on the domain $D = \{(x_1, \ldots, x_n) | 0<x_i<\pi, \sum_1^n x_i=P   \}.$ Then $f(x)\le f(x^*)$ for all $x\in D$ where   
$x^* = (\frac{P}{n}, \ldots, \frac{P}{n}).$  
  
To prove this claim we will show that starting from any point $x\in D , x\ne x^*$ we can build a sequence $\{x^k\}_{k=1..\infty}, x^k\in D,$ such that:

i) $x^1 = x.$  

ii) $f(x)$ is non-decreasing on the segment connecting $x^i$ and $x^{i+1}.$ 

iii) $\lim_{k\rightarrow \infty} x^k =x^*.$  

Once we have done this, it follows immediately that $ f(x) \le lim_{k\rightarrow \infty} f(x^k) = f(x^*).$ 

We construct the sequence $\{x^k\}$ as follows. Put $x^1=x.$ Suppose that we already have $x^k,$ we will find $x^{k+1}$ as follows:

- If $x^k=x^*,$ put $x^{k+1} = x^*.$

- If $x^k = (x^k_1, \ldots, x^k_n)  \ne x^*,$ suppose that  $x^k_{m} = \min_{i=1..n} x^k_i\  \textrm{and} \ x^k_{M} = \max_{i=1..n} x^k_i,$ we then specified $x^{k+1}$ by :

$$ x^{k+1}_i := {\Bigg \{ }  \begin{matrix}  x^k_i  & \textrm{for}& \  i\ne m, M\\  \frac{x^k_m+x^k_M}{2}  & \textrm{for}& \  i = m\  \textrm{or}\   i = M. \end{matrix}$$       
   
It is not hard to see that  the sequence constructed this way converges to $x^*.$ It remains to prove that it satisfies  (ii) above. Without loss of generality, it is enough to show that if $x=(x_1,x_2, \ldots, x_n)$ and $x_1 =  \max_{i=1..n} x_i,$ 
$x_2 =  \min_{i=1..n} x_i$ then $f(x)$ is non-decreasing on the segment $l(t):= (x_1-\frac{t(x_1-x_2)}{2}, x_2+\frac{t(x_1-x_2)}{2}, x_3, \ldots, x_n), \ 0\le t\le 1.$

When restricted to the segment $l(t),$ our function has the form :
\begin{eqnarray*}
f(t) &=& \sum_{k=1}^\infty \frac{\sin k(x_1-\frac{t(x_1-x_2)}{2})\sin k(x_2+\frac{t(x_1-x_2)}{2}) \ldots \sin kx_n}{k^{n-2}}\\
& = & \sum_{k=1}^\infty \frac{(\cos k(x_1-x_2 - t(x_1-x_2))- \cos k(x_1+x_2))\sin kx_3  \ldots \sin kx_n}{2k^{n-2}}.
\end{eqnarray*}
We compute :
$$f'(t) = (x_1-x_2)\sum_{k=1}^\infty \frac{\sin k(x_1-x_2 - t(x_1-x_2)) \sin kx_3  \ldots \sin kx_n}{2k^{n-3}}.$$

Since $x_1>x_2,$ by Corollary 2.2, we deduce that  $f'(t)\ge 0,\ \forall 0\le t \le 1.$ So we obtain the first part of the theorem.

For $n$ even, by Corollary 2.2, it is easy to see that $f'(t)$ above is strictly positive for $x$ sufficiently close to $x^*.$ We conclude that $x^*$ is the unique point where $f$ attains its maximum value.

For $n$ odd, in order for the function $f(x)$  not to be identically zero on the domain $D,$ by Corollary 2.2, we must have $P \le (n-1)\pi.$ Even with this assumption, in certain cases, $x^*$ may not be  the unique point where $f$ attains its maximum value. For example, when $n=5$ and $P=4\pi,$ then it is not hard to check that 
$$f(\frac{4\pi}{5},\frac{4\pi}{5},\frac{4\pi}{5},\frac{4\pi}{5},\frac{4\pi}{5}) = f(\frac{4\pi}{5}+\epsilon,\frac{4\pi}{5}-\epsilon,\frac{4\pi}{5},\frac{4\pi}{5},\frac{4\pi}{5})$$  
for all $\epsilon >0$ sufficiently small.

\textit {Proof of Theorem 1.2.}
It follows from the proof of Theorem 1.1 that a regular polygon  is the unique one 
which is the most flexible in the set of all polygons with the same perimeter.    
So, by Theorem 1.1, it is enough to restrict ourself to the regular polygons.
 
Let $0<x<\pi,$ and denote by $V_n(x)$ the symplectic volume of the moduli space of regular spherical $n$-gons with side $x$. From Theorem 2.1, we know that 
$$ V_n(x)= \frac{2^{n-1}}{\pi} \sum_{k=1}^\infty \frac{(\sin kx)^n}{k^{n-2}}.$$ 
As
 $$V_n'(x)= \frac{2^{n-1}}{\pi} \sum_{k=1}^\infty  \frac{n(\sin kx)^{n-1}\cos kx}{k^{n-3}}= \frac{2^{n-2}n}{\pi} \sum_{k=1}^\infty  \frac{(\sin kx)^{n-2}\sin 2kx}{k^{n-3}}.$$ 
If $0 <  x < \pi/2,$  Corollary 2.2 tells us that $V'_n(x) \ge 0.$ 

In the case  $\pi/2 <  x < \pi,$  we write :
$$V_n'(x) =  \frac{2^{n-2}n}{\pi} \sum_{k=1}^\infty - \frac{(\sin kx)^{n-2}\sin 2k(\pi-x)}{k^{n-3}}$$
 and again by Corollary 2.2, we get $V_n'(x) \le 0.$   

Notice that $V_n'(\pi/2)=0$ and moreover, by using the second part of Corollary 2.2, it is not hard to check that $V_n'(x)$ 
is strictly positive (resp. negative) when $x$ is sufficiently close to $\pi/2$ on the left (resp. right). From all this we conclude that $\pi/2$ is the unique global maximum of $V_n(x)$ on the interval $(0, \pi).$ 
 So,  we get the conclusion of  Theorem 1.2. 
\begin{remark}
Theorem 1.2 seems to agree with the intuitive reasoning that as the side-lengths are $\pi/2,$ each vertex moves on the great sphere, i.e., on the biggest possible space, therefore the polygon is the most flexible.  
\end{remark}

\end{document}